\documentclass[a4paper,12pt]{article}
\usepackage{amsmath}
\usepackage{color}
\usepackage{amssymb}
\usepackage{latexsym}
\numberwithin{equation}{section}

\def\R{{\bf R}}

\def\N{{\bf N}}

\def\d{\displaystyle}
\def\e{{\varepsilon}}

\def\wt{\widetilde}

\def\v#1{\mbox{\boldmath $#1$}}

\newtheorem{thm}{Theorem}[section]

\newtheorem{lem}{Lemma}[section]

\newtheorem{rem}{Remark}[section]
\newtheorem{Def}{Definition}[section]

\title{Blow-up for semilinear wave equations\\
with the scale invariant damping \\ and super-Fujita exponent
}
\author{
Ning-An Lai
\footnote{Department of Mathematics,
Lishui University,
Lishui City 323000,
China.
e-mail: hyayue@gmail.com.}
\quad
Hiroyuki Takamura
\footnote{Department of Complex and Intelligent Systems,
Faculty of Systems Information Science,
Future University Hakodate,
116-2 Kamedanakano-cho,
Hakodate, Hokkaido 041-8655, Japan.
e-mail: takamura@fun.ac.jp.}
\quad
Kyouhei Wakasa
\footnote{
College of Liberal Arts, Mathematical Science Research Unit,
Muroran Institute of Technology, 27-1, Mizumoto-cho,
Muroran, Hokkaido 050-8585, Japan.
email: wakasa@mmm.muroran-it.ac.jp.
}
}
\date{
\[
\begin{array}{ll}
\mbox{\footnotesize{\bf Keywords:}}
& \mbox{\footnotesize damped wave equation, semilinear, blow-up}\\
\mbox{\footnotesize{\bf MSC2010:}}
& \mbox{\footnotesize primary 35L71, secondary 35B44}\\
\end{array}
\]
}

\begin{document}
\maketitle
\begin{abstract}
The blow-up for semilinear wave equations with the scale
invariant damping has been well-studied for sub-Fujita exponent.
However, for super-Fujita exponent, there is only one blow-up result
which is obtained  in 2014 by Wakasugi in the case of non-effective damping.
In this paper we extend his result in two aspects by showing that:
(I) the blow-up will happen for bigger exponent, which is closely related to the Strauss exponent, the critical number for non-damped semilinear wave equations;
(II) such a blow-up result is established for a wider range of the constant
than the known non-effective one  in the damping term.
\end{abstract}


\section{Introduction}
\par\quad
In this paper, we consider the following initial value problem.
\begin{equation}
\label{IVP}
\left\{
\begin{array}{l}
\d u_{tt}-\Delta u+\frac{\mu}{1+t}u_t=|u|^p
\quad \mbox{in}\ \R^n\times[0,\infty),\\
u(x,0)=\e f(x),\ u_t(x,0)=\e g(x), \quad  x\in\R^n,
\end{array}
\right.
\end{equation}
where $\mu>0,\ f,g\in C_0^{\infty}(\R^n)$ and $n\in\N$.
We assume that $\e>0$ is a \lq\lq small" parameter.
\par
First, we shall outline a background of (\ref{IVP}) briefly
according to the classifications by Wirth in \cite{Wir1, Wir2, Wir3}
for the corresponding linear problem.
Let $u^0$ be a solution of the initial value problem
for the following linear damped wave equation.
\begin{equation}
\label{scale_eq}
\left\{
\begin{array}{l}
\d u^0_{tt}-\Delta u^0+\frac{\mu}{(1+t)^\beta}u^0_t=0,
\quad \mbox{in}\ \R^n \times[0,\infty),\\
u^0(x,0)=u_1(x),\ u^0_t(x,0)=u_2(x), \quad x\in\R^n,
\end{array}
\right.
\end{equation}
where $\mu>0$, $\beta\in\R$, $n\in\N$ and $u_1,u_2\in C^\infty_0(\R^n)$.
When $\beta\in(-\infty,-1)$, we say that the damping term is \lq\lq overdamping"
in which case the solution does not decay to zero when $t\rightarrow\infty$.
When $\beta\in[-1,1)$, the solution behaves like that of the heat equation,
which means that the term $u^0_{tt}$ in (\ref{scale_eq}) has no influence
on the behavior of the solution.
In fact, $L^p$-$L^q$ decay estimates of the solution
which are almost the same as those of the heat equation are established.
In this case, we say that the damping term is \lq\lq effective."
In contrast, when $\beta\in(1,\infty)$,
it is known that the solution behaves like that of the wave equation,
which means that the damping term in (\ref{scale_eq}) has no influence
on the behavior of the solution.
In fact, in this case the solution scatters to that of the free wave equation
when $t\rightarrow\infty$,
and thus we say that we have \lq\lq scattering."
When $\beta=1$, the equation in (\ref{scale_eq}) is invariant under the following scaling
\[
\wt{u^0}(x,t):=u^0(\sigma x, \sigma(1+t)-1),\ \sigma>0,
\]
and hence we say that the damping term is \lq\lq scale invariant."
The remarkable fact in this case is that the behavior of the solution of (\ref{scale_eq})
is determined by the value of $\mu$.
Actually, for $\mu\in(0,1)$, it is known that the asymptotic behavior of the solution
is closely related to that of the free wave equation.
For this range of $\mu$,
we say that the damping term is \lq\lq non-effective."
However, the threshold of $\mu$ according to the behavior of the solution is still open.
We conjecture that it may be $\mu=1$ since we have the following $L^2$ estimates:
\[
\|u^0(\cdot,t)\|_{L^2}\lesssim
\|u_1\|_{L^2}+\|u_2\|_{H^{-1}}\times
\left\{
\begin{array}{ll}
(1+t)^{1-\mu}& \mbox{if}\ \mu \in(0,1),\\
\log (e+t) &  \mbox{if}\ \mu=1,\\
1 & \mbox{if}\ \mu>1.
\end{array}
\right.
\]
In this way, we may summarize all the classifications of the damping term in (\ref{scale_eq}) 
in the following table.
\begin{center}
\begin{tabular}{|c|c|}
\hline
$\beta\in(-\infty,-1)$ & overdamping\\
\hline
$\beta\in[-1,1)$ & effective\\
\hline
$\beta=1$ &
\begin{tabular}{c}
scaling invariant\\
$\mu\in(0,1)\Rightarrow$ non-effective
\end{tabular}
\\
 \hline
$\beta\in(1,\infty)$ & scattering\\
\hline
\end{tabular}
\end{center}
\par
Next, we consider the following initial value problem for semilinear damped wave equation.
\begin{equation}
\label{scale_eq_nonlinear}
\left\{
\begin{array}{l}
\d v_{tt}-\Delta v+\frac{\mu}{(1+t)^\beta}v_t=|v|^p, \quad \mbox{in}\ \R^n \times[0,\infty),\\
v(x,0)=\e f(x),\ v_t(x,0)=\e g(x), \quad x\in\R^n,
\end{array}
\right.
\end{equation}
where $\mu>0,\ \beta\ge-1,\ f,g\in C_0^{\infty}(\R^n)$ and $n\in\N$.
We assume that $\e>0$ is a \lq\lq small" parameter.
\par
For the constant coefficient case, $\beta=0$,
Todorova and Yordanov \cite{TY} have shown that
the energy solution of (\ref{scale_eq_nonlinear}) exists globally-in-time
for \lq\lq small" initial data if $p>p_F(n)$, where
\begin{equation}
\label{fujita}
 p_F(n):=1+\frac{2}{n}
 \end{equation}
is the so-called Fujita exponent, the critical exponent for semilinear heat equations.
It has been also obtained in \cite{TY} that
the solution of (\ref{scale_eq_nonlinear})
blows-up in finite time for some positive data if $1<p<p_F(n)$.
The critical case $p=p_F(n)$ has been studied by Zhang \cite{Z}
by showing the blow-up result.
We note that Li and Zhou \cite{LZ}, or Nishihara \cite{N},
have obtained the sharp upper bound of the lifespan
which is the maximal existence time of solutions of (\ref{scale_eq_nonlinear})
in the case of $n=1,2$,  or $n=3$, respectively.
The sharpness of the upper bound has been studied by Li \cite{Li}
including the result for more general equations with all $n\ge1$,
but for smooth nonlinear terms.
The sharp lower bound has been obtained by Ikeda and Ogawa \cite{IO} for the critical case.
Recently,  Lai and Zhou \cite{LaZ} have obtained the sharp upper bound of the lifespan
in the critical case for $n\ge4$.
\par
For the variable coefficient case of the most part of the effective damping
with $-1<\beta<1$, Lin, Nishihara and Zhai \cite{LNZ12} have obtained
the blow-up result if $1<p\le p_F(n)$
and the global existence result if $p>p_F(n)$.
Later, D'Abbicco, Lucente and Reissig \cite{DLR13} have extended
the global existence result to more general equations.
For the precise estimates of the lifespan in this case,
see Introduction in Ikeda and Wakasugi \cite{IW}.
Recently, similar results on the remaining part of effective damping with $\beta=1$
have been obtained by Wakasugi \cite{WY17} for the global existence part,
and by Fujiwara, Ikeda and  Wakasugi \cite{FIW} for the blow-up part.
The sharp estimates of the lifespan are also obtained by \cite{FIW}
except for the upper bound in the critical case.
\par
Now, let us turn back to our problem (\ref{IVP}).
Wakasugi \cite{WY14_scale} has obtained the blow-up result
if $1<p\le p_F(n)$ and $\mu>1$, or $1<p\le 1+2/(n+\mu-1)$ and $0<\mu\le 1$.
He has also shown in \cite{WY14} that an upper bound of the lifespan is
\[
\left\{
\begin{array}{ll}
C\e^{-(p-1)/\{2-n(p-1)\}} & \mbox{if}\ 1<p<p_F(n)\ \mbox{and}\ \mu\ge1,\\
C\e^{-(p-1)/\{2-(n+\mu-1)(p-1)\}} & \mbox{if}\  \d1<p<1+\frac{2}{n+\mu-1}
\ \mbox{and}\ 0<\mu<1,
\end{array}
\right.
\]
where $C$ is a positive constant independent of $\e$.
We note that the both proofs in \cite{WY14} and \cite{WY14_scale}
are based on the so-called \lq\lq test function method" introduced by Zhang \cite{Z}.
On the other hand, D'Abbicco \cite{DABI} has obtained
the global existence result if $p>p_F(n)$,
but $\mu$ has to satisfy
\[
\mu\ge
\left\{
\begin{array}{cl}
 5/3 & \mbox{for}\  n=1,\\
 3 & \mbox{for}\ n=2,\\
n+2 & \mbox{for}\  n\ge3\ \left(\mbox{and}\ p\le1+2/(n-2)\right).
\end{array}
\right.
\]
It is remarkable that, by the so-called Liouville transform;
\[
w(x,t):=(1+t)^{\mu/2}u(x,t),
\]
(\ref{IVP}) can be rewritten as
\begin{equation}
\label{IVP_1}
\left\{
\begin{array}{l}
\d w_{tt}-\Delta w+\frac{\mu(2-\mu)}{4(1+t)^2}w=\frac{|w|^p}{(1+t)^{\mu(p-1)/2}}
\quad \mbox{in}\ \R^n \times[0,\infty),\\
w(x,0)=\e f(x),\ w_t(x,0)=\e \{(\mu/2)f(x)+g(x)\}, \quad x\in\R^n.
\end{array}
\right.
\end{equation}
When $\mu=2$, D'Abbicco, Lucente and Reissig \cite{DLR14} have obtained
the following result.
Let
\begin{equation}
\label{critical_value_IVP1}
p_c(n):=\max\left\{p_{F}(n),\ p_0(n+2)\right\},
\end{equation}
where
\begin{equation}
\label{strauss}
p_0(n):=\frac{n+1+\sqrt{n^2+10n-7}}{2(n-1)}
\quad(n\ge2)
\end{equation}
is the so-called Strauss exponent, the positive root of the quadratic equation,
\begin{equation}
\label{gamma}
\gamma(p,n):=2+(n+1)p-(n-1)p^2=0.
\end{equation}
We note that $p_0(n)$ is the critical exponent for semilinear wave equations,
$\mu=0$ in (\ref{IVP}).
They have shown in \cite{DLR14} that the problem (\ref{IVP}) admits a global-in-time solution
in the classical sense for \lq\lq small" $\e$ if $p>p_c(n)$ in the case of $n=2,3$
although the radial symmetry is assumed in $n=3$,
and that the classical solution of (\ref{IVP}) with positive data
blows-up in finite time if $1<p\le p_c(n)$ and $n\ge1$. 
In the same year, with radial symmetric assumption,
D'Abbicco and Lucente \cite{D-L}
extended the global existence result for $p_0(n+2)<p<1+2/(\max\{2,(n-3)/2\})$
to odd higher dimensions ($n\geq 5$).
We remark that, in the case of $n=1$,
Wakasa \cite{wak16} has studied the estimates of the lifespan
and has shown that the critical exponent $p_c(1)=p_F(1)=3$
changes to $p_0(1+2)=1+\sqrt{2}$
when the nonlinearity is a sign-changing type, $|u|^{p-1}u$,
and the initial data is of odd functions.
Both results in \cite{DLR14} and \cite{wak16} heavily rely
on the special structure of the massless wave equations,
$\mu=2$ in (\ref{IVP_1}).
In view of them, $\mu=2$ may be an exceptional case.
Recalling Wirth' classification in the linear problem, (\ref{scale_eq}),
one may regard $\mu=1$ as a threshold also for the semilinear problem, (\ref{IVP}).
In this sense, the blow-up result in Wakasugi \cite{WY14_scale} says that
the solution may be \lq\lq heat-like" if $\mu>1$.
Here, \lq\lq heat-like" means that the critical exponent for (\ref{IVP}) is Fujita exponent.
\par
In this paper, we claim that the solution of (\ref{IVP}) is
\lq\lq wave-like" in some case even for $\mu>1$.
Here, \lq\lq wave-like" means that
the critical exponent for (\ref{IVP}) is bigger than Fujita exponent
and is related to Strauss exponent.
We also conjecture that such a threshold of $\mu$ depends on the space dimension $n$.
The main tool of our result is Kato's lemma in Kato \cite{Kato80}
on ordinary differential inequalities
which is improved to be applied to semilinear wave equations
by Takamura \cite{Takamura15}.
Together with Yordanov and Zhang's estimate in \cite{YZ06},
we can prove a new blow-up result for wave-like solutions
by means of some special transform for the time-derivative
of the spatial integral of unknown functions.
\par
This paper is organized as follows.
In the next section, we state our main result.
In the section 3, we estimate the spatial integral of unknown functions from below.
Making use of such an estimate, we prove the main result for $\mu\ge2$ in section 4,
and for $0<\mu<2$ in section 5.

\section{Main Result}
\label{section:lifespan}
\par\quad
First we define an energy solution of (\ref{IVP}).
\begin{Def}
We say that $u$ is an energy solution of (\ref{IVP}) on $[0,T)$
if
\begin{equation}
\label{reg}
u\in C([0,T),H^1(\R^n))\cap C^1([0,T),L^2(\R^n))\cap L_{\rm loc}
^p(\R^n\times[0,T))
\end{equation}
satisfies
\begin{equation}
\label{sol}
\begin{array}{l}
\d\int_{\R^n}u_t(x,t)\phi(x,t)dx-\int_{\R^n}u_t(x,0)\phi(x,0)dx\\
\d+\int_0^tds\int_{\R^n}\left\{-u_t(x,s)\phi_t(x,s)+\nabla u(x,s)\cdot\nabla\phi(x,s)\right\}dx\\
\d+\int_0^tds\int_{\R^n}\frac{\mu u_t(x,s)}{1+s}\phi(x,s)dx
=\int_0^tds\int_{\R^n}|u(x,s)|^p\phi(x,s)dx
\end{array}
\end{equation}
with any $\phi\in C_0^{\infty}(\R^n\times[0,T))$ and any $t\in[0,T)$.
\end{Def}

We note that, employing the integration by parts in (\ref{sol})
and letting $t\rightarrow T$, we have that
\[
\begin{array}{l}
\d\int_{\R^n\times[0,T)}
u(x,s)\left\{\phi_{tt}(x,s)-\Delta \phi(x,s)
-\left(\frac{\mu\phi(x,s)}{1+s}\right)_s \right\}dxds\\
\d=\int_{\R^n}\mu u(x,0)\phi(x,0)dx-\int_{\R^n}u(x,0)\phi_t(x,0)dx\\
\d\quad+\int_{\R^n}u_t(x,0)\phi(x,0)dx+\int_{\R^n\times[0,T)}|u(x,s)|^p\phi(x,s)dxds.
\end{array}
\]
This is exactly the definition of a weak solution of (\ref{IVP}).
\par
Our main result is the following theorem.

\begin{thm}
\label{thm:main}
Let $n\ge2$,
\begin{equation}
\label{mup}
0<\mu<\mu_0(n):=\frac{n^2+n+2}{2(n+2)}
\quad\mbox{and}\quad
p_F(n)\le p<p_0(n+2\mu).
\end{equation}
Assume that both $f\in H^1(\R^n)$ and $g\in L^2(\R^n)$
are non-negative and do not vanish identically.
Suppose that an energy solution $u$ of (\ref{IVP}) satisfies
\begin{equation}
\label{support}
\mbox{\rm supp}\ u\ \subset\{(x,t)\in\R^n\times[0,T)\ :\ |x|\le t+R\}
\end{equation}
with some $R\ge1$.
Then, there exists a constant $\e_0
=\e_0(f,g,n,p,\mu,R)>0$ such that $T$ has to satisfy
\begin{equation}
\label{lifespan}
T\le C\e^{-2p(p-1)/\gamma(p,n+2\mu)}
\end{equation}
for $0<\e\le\e_0$, where $C$ is a positive constant independent of $\e$.
\end{thm}
\begin{rem}
Theorem \ref{thm:main} can be established also for $n=1$
if one define $\phi_1(x)=e^x+e^{-x}$ for $x\in\R$ in Section 3 below.
But the result is not new.
See the following two remarks.
\end{rem}
\begin{rem}
In view of (\ref{fujita}) and (\ref{gamma}),
one can see that
\[
\gamma\left(p_F(n),n+2\mu\right)
=\frac{2}{n^2}\left\{n^2+n+2-2(n+2)\mu\right\}.
\]
Therefore $0<\mu<\mu_0(n)$ is equivalent to
\[
p_F(n)<p_0(n+2\mu).
\]
We note that $\mu_0(2)=1$.
This means that Theorem \ref{thm:main} just covers
the non-effective range of $\mu$ for $n=2$.
Since $\mu_0(n)$ is increasing in $n$,
Theorem \ref{thm:main} gives us the blow-up result on super-Fujita exponent
even for $\mu$ in outside of the non-effective range for $n\ge3$.
We also note that
$\mu_0(n)<2$ for $n=2,3,4$ and
$\mu_0(n)>2$ for $n\ge5$.
\end{rem}

\begin{rem}
One can see also that
\[
\gamma\left(1+\frac{2}{n+\mu-1},n+2\mu\right)\
=\frac{2\left\{(n-1)^2+(n-3)\mu\right\}}{(n+\mu-1)^2}.
\]
Therefore we have that
\[
1+\frac{2}{n+\mu-1}<p_0(n+2\mu)
\]
for $n=2$ and $0<\mu<1$, or $n\ge3$ and $\mu>0$.
This means that Theorem \ref{thm:main}
includes the blow-up result in Wakasugi \cite{WY14_scale}.
\end{rem}

\begin{rem}
If $\beta$ is in the scattering range, $(1,\infty)$, for the problem,
\[
u_{tt}-\Delta u+\frac{\mu}{(1+t)^\beta}u_t=|u|^p,
\]
the result will be
\[
T\le C\e^{-2p(p-1)/\gamma(p,n)}
\quad\mbox{for}\ 1<p<p_0(n)
\]
for all $\mu>0$.
This estimate coincides with the one for non-damped equation,
\[
u_{tt}-\Delta u=|u|^p
\]
except for the case of $\int_{\R^2}g(x)dx\neq0$ in $n=2$ and $1<p\le2$.
See Introduction of Takamura \cite{Takamura15} for its summary.
The proof of this fact will appear in our forthcoming paper.
\end{rem}
\section{Lower bound of the functional}
\par\quad
Let $u$ be an energy solution of (\ref{IVP}) on $[0,T)$.
We estimate
\[
F_0(t):=\int_{\R^n}u(x,t)dx
\]
from below in this section.
Choosing the test function $\phi=\phi(x,s)$ in (\ref{sol}) to satisfy
$\phi\equiv 1$ in $\{(x,s)\in \R^n\times[0,t]:|x|\le s+R\}$, we get
\[
\begin{array}{l}
\d\int_{\R^n}u_t(x,t)dx-\int_{\R^n}u_t(x,0)dx\\
\d+\int_0^tds\int_{\R^n}\frac{\mu u_t(x,s)}{1+s}dx
=\int_0^tds\int_{\R^n}|u(x,s)|^pdx,
\end{array}
\]
which means that
\[
F'_0(t)-F'_0(0)+\int_0^t\frac{\mu F'_0(s)}{1+s}dx
=\int_0^tds\int_{\R^n}|u(x,s)|^pdx.
\]
All the quantities in this equation except for $F'_0(t)$ are differentiable in $t$,
so that so is $F'_0(t)$.
Hence we have
\begin{equation}
\label{ode}
F''_0(t)+\frac{\mu F'_0(t)}{1+t}
=\int_{\R^n}|u(x,t)|^pdx.
\end{equation}
Integrating this equation with a multiplication by $(1+t)^\mu$, we obtain
\begin{equation}
\label{F'_0}
(1+t)^\mu F'_0(t)-F'_0(0)=\int_0^t(1+s)^{\mu}ds\int_{\R^n}|u(x,s)|^pdx.
\end{equation}
It follows from this equation and the assumption on the initial data that
\begin{equation}
\label{positive}
F'_0(t)\ge (1+t)^{-\mu}F'_0(0)>0
\quad\mbox{and}\quad
F_0(t)\ge F_0(0)>0
\quad\mbox{for}\quad
t\ge0.
\end{equation}
\par
From now on, we employ the modified argument of Yordanov and Zhang \cite{YZ06}.
Let us define
\[
F_1(t):=\int_{\R^n}u(x,t)\psi_1(x,t)dx,
\]
where
\[
\psi_1(x,t):=\phi_1(x)e^{-t},
\quad
\phi_1(x):=\int_{S^{n-1}}e^{x\cdot\omega}dS_\omega.
\]
In view of (\ref{F'_0}) and the argument of (2.4)-(2.5) in \cite{YZ06},
we know that there is a positive constant $C_1=C_1(n,p,R)$ such that
\begin{equation}
\label{YZ}
(1+t)^\mu F'_0(t)-F'_0(0)\ge C_1\int_0^t(1+s)^{\mu+(n-1)(1-p/2)}|F_1(s)|^pds.
\end{equation}
\par
In order to get a lower bound of $F_1(t)$,
we turn back to (\ref{sol}) and obtain that
\[
\begin{array}{l}
\d\frac{d}{dt}\int_{\R^n}u_t(x,t)\phi(x,t)dx\\
\d+\int_{\R^n}\left\{-u_t(x,t)\phi_t(x,t)-u(x,t)\Delta\phi(x,t)\right\}dx\\
\d+\int_{\R^n}\frac{\mu u_t(x,t)}{1+t}\phi(x,t)dx
=\int_{\R^n}|u(x,t)|^p\phi(x,t)dx.
\end{array}
\]
Multiplying the above equality by$(1+t)^{\mu}$, we have that
\[
\begin{array}{l}
\d\frac{d}{dt}\left\{(1+t)^\mu\int_{\R^n}u_t(x,t)\phi(x,t)dx\right\}\\
\d+(1+t)^\mu\int_{\R^n}\left\{-u_t(x,t)\phi_t(x,t)-u(x,t)\Delta\phi(x,t)\right\}dx\\
\d=(1+t)^\mu\int_{\R^n}|u(x,t)|^p\phi(x,t)dx.
\end{array}
\]
Integrating this equality over $[0,t]$, we get
\[
\begin{array}{l}
\d(1+t)^\mu\int_{\R^n}u_t(x,t)\phi(x,t)dx-\e\int_{\R^n}g(x)\phi(x,0)dx\\
\d-\int_0^tds\int_{\R^n}(1+s)^\mu u_t(x,s)\phi_t(x,s)dx\\
\d=\int_0^tds\int_{\R^n}\left\{(1+s)^\mu u(x,s)\Delta\phi(x,s)+(1+s)^\mu|u(x,s)|^p\phi(x,s)\right\}dx.
\end{array}
\]
It follows from this equation and
\[
\begin{array}{l}
\d\int_0^tds\int_{\R^n}(1+s)^\mu u_t(x,s)\phi_t(x,s)dx\\
\d=(1+t)^\mu\int_{\R^n}u(x,t)\phi_t(x,t)dx-\int_{\R^n}u(x,0)\phi_t(x,0)dx\\
\quad\d-\int_0^tds\int_{\R^n}\mu(1+s)^{\mu-1}u(x,s)\phi_t(x,s)dx\\
\quad\d-\int_0^tds\int_{\R^n}(1+s)^\mu u(x,s)\phi_{tt}(x,s)dx,
\end{array}
\]
which follows from integration by parts that
\[
\begin{array}{l}
\d(1+t)^\mu\int_{\R^n}\left\{u_t(x,t)\phi(x,t)-u(x,t)\phi_t(x,t)\right\}dx\\
\d-\e\int_{\R^n}g(x)\phi(x,0)dx+\e\int_{\R^n}f(x)\phi_t(x,0)dx\\
\d+\int_0^tds\int_{\R^n}\mu(1+s)^{\mu-1}u(x,s)\phi_t(x,s)dx\\
\d=\int_0^tds\int_{\R^n}(1+s)^\mu u(x,s)\{\Delta\phi(x,s)-\phi_{tt}(x,s)\}dx\\
\quad\d+\int_0^tds\int_{\R^n}(1+s)^\mu|u(x,s)|^p\phi(x,s)dx.
\end{array}
\]
\par
If we put
\[
\phi(x,t)=\psi_1(x,t)=e^{-t}\phi_1(x)
\quad\mbox{on}\quad \mbox{supp}\ u,
\]
we have
\[
\phi_t=-\phi,\ \phi_{tt}=\Delta\phi \quad\mbox{on}\quad \mbox{supp}\ u.
\]
Hence we obtain that
\[
\begin{array}{l}
\d(1+t)^\mu F_1'(t)+2(1+t)^\mu F_1(t)-\e\int_{\R^n}\left\{f(x)+g(x)\right\}\phi(x)dx\\
\d=\int_0^t\mu(1+s)^{\mu-1}F_1(s)ds
+\int_0^tds\int_{\R^n}(1+s)^\mu|u(x,s)|^p\phi(t, x)dx,
\end{array}
\]
which yields
\[
F_1'(t)+2F_1(t)\ge\frac{C_{f,g}\e}{(1+t)^\mu}+\frac{1}{(1+t)^\mu}\int_0^t\mu(1+s)^{\mu-1}F_1(s)ds,
\]
where
\[
C_{f,g}:=\int_{\R^n}\left\{f(x)+g(x)\right\}\phi_1(x)dx>0.
\]
Integrating this inequality over $[0,t]$ with a multiplication by $e^{2t}$, we get
\begin{equation}
\label{ineq:F1}
\begin{array}{ll}
e^{2t}F_1(t)\ge
&\d F_1(0)+C_{f,g}\e\int_0^t\frac{e^{2s}}{(1+s)^\mu}ds\\
&\d+\int_0^t\frac{e^{2s}}{(1+s)^\mu}ds\int_0^s\mu(1+r)^{\mu-1}F_1(r)dr.
\end{array}
\end{equation}
We note that the assumption on $f$ implies $F_1(0)>0$.
Hence we find that there is no zero point of $F_1(t)$ for $t>0$.
Because the continuity of $F_1$ implies $F_1(t)>0$ for small $t>0$.
If one assumes that there is a nearest zero point $t_0$ of $F_1$ to $0$,
then one has a contradiction in (\ref{ineq:F1});
\[
\begin{array}{ll}
e^{2t_0}F_1(t_0)=0\ge
&\d F_1(0)+C_{f,g}\e\int_0^{t_0}\frac{e^{2s}}{(1+s)^\mu}ds\\
&\d+\int_0^{t_0}\frac{e^{2s}}{(1+s)^\mu}ds\int_0^s\mu(1+r)^{\mu-1}F_1(r)dr.
\end{array}
\]
The last term in the right-hand side of this inequality is positive by $F_1(t)>0$ for $0<t<t_0$.
Turning back to (\ref{ineq:F1}), we obtain
\[
F_1(t)>\frac{C_{f,g}\e}{2(1+t)^\mu}(1-e^{-2t})+e^{-2t}F_1(0)
\ge\frac{C_{f,0}\e}{2(1+t)^\mu}\quad\mbox{for}\ t\ge0.
\]
Here we have used the fact that $C_{f,g}>C_{f,0}$.
\par
Plugging this estimate into (\ref{YZ}), we have
\[
(1+t)^\mu F'_0(t)-F'_0(0)>\frac{C_1C_{f,0}^p}{2^p}\e^p\int_0^t(1+s)^{\mu(1-p)+(n-1)(1-p/2)}ds.
\]
Since $F'_0(0)>0$ and it follows from $t\ge1$ that
\[
\begin{array}{l}
\d\int_0^t(1+s)^{\mu(1-p)+(n-1)(1-p/2)}ds
\ge(2t)^{\mu(1-p)}\int_{t/2}^ts^{(n-1)(1-p/2)}ds\\
\qquad\d\ge2^{\mu(1-p)-(n-1)(1-p/2)+1}t^{1+\mu(1-p)+(n-1)(1-p/2)},
\end{array}
\]
we obtain that
\[
F'_0(t)>C_2\e^pt^{1-\mu p+(n-1)(1-p/2)}
\quad\mbox{for}\ t\ge1,
\]
where
\[
C_2:=\frac{C_1C_{f,0}^p}{2^{(\mu+1)p+(n-1)(1-p/2)-1}}>0.
\]
Here we have used the fact that
\[
1+(n-1)\left(1-\frac{p}{2}\right)>0
\]
follows from
\[
p<p_0(n+2\mu)<
\left\{
\begin{array}{ll}
p_0(n)\le p_0(4)=2 & \mbox{for}\ n\ge4,\\
p_0(3)=1+\sqrt{2}<3 & \mbox{for}\ n=3,\\
p_0(2)=(3+\sqrt{17})/2<4 & \mbox{for}\ n=2.
\end{array}
\right.
\]
Integrating this inequality over $[1,t]$ and making use of $F_0(0)>0$ and
\[
\int_1^ts^{1-\mu p+(n-1)(1-p/2)}ds
\ge t^{-\mu p}\int_{t/2}^ts^{1+(n-1)(1-p/2)}ds
\quad\mbox{for}\ t\ge2,
\]
we get
\begin{equation}
\label{lower bound}
F_0(t)>C_3\e^pt^{2-\mu p+(n-1)(1-p/2)}
\quad\mbox{for}\ t\ge2,
\end{equation}
where
\[
C_3:=\frac{C_2}{2^{2+(n-1)(1-p/2)}}>0.
\]


\section{Proof of Theorem \ref{thm:main} for $\v{\mu\ge2}$}
\par\quad
Let us define
\[
F(t):=\int_{\R^n}w(x,t)dx=(1+t)^{\mu/2}F_0(t),
\]
where $w$ is the solution of (\ref{IVP_1}).
We note that (\ref{ode}) yields
\begin{equation}
\label{ineq_LG}
F''(t)+\frac{\mu(2-\mu)}{4(1+t)^2}F(t)
=(1+t)^{-\mu(p-1)/2}\int_{\R^n}|w(x,t)|^pdx.
\end{equation}
Then it follows from (\ref{support}) and H\"older's inequality that
\begin{equation}
\label{ineq_LF}
\int_{\R^n}|w(x,t)|^pdx\ge\{\mbox{vol}({\bf B}^n(0,1))\}^{1-p}(t+R)^{-n(p-1)}|F(t)|^p.
\end{equation}
By combining \eqref{ineq_LG} and \eqref{ineq_LF}, and noting the assumption $R\ge1$, we come to
\begin{equation}
\label{ineq_F}
F''(t)+\frac{\mu(2-\mu)}{4(1+t)^2}F(t)
\ge C_4(1+t)^{-(n+\mu/2)(p-1)}|F(t)|^p
\quad\mbox{for}\ t\ge0,
\end{equation}
where
\[
C_4:=\{\mbox{vol}({\bf B}^n(0,1))\}^{1-p}R^{-n(p-1)}>0.
\]
Due to (\ref{positive}), we have that
\begin{equation}
\label{positive_F}
\begin{array}{l}
F(t)=(1+t)^{\mu/2}F_0(t)>0,\\
\d F'(t)=\frac{\mu}{2}(1+t)^{\mu/2-1}F_0(t)+(1+t)^{\mu/2}F'_0(t)>0,
\end{array}
\end{equation}
which implies that
\begin{equation}
\label{data_Fat0}
\begin{array}{l}
F(0)=F_0(0)=\|f\|_{L^1(\R^n)}\e,\\
\d F'(0)=\frac{\mu}{2}F_0(0)+F'_0(0)
=\left(\frac{\mu}{2}\|f\|_{L^1(\R^n)}+\|g\|_{L^1(\R^n)}\right)\e.
\end{array}
\end{equation}
\par
From now on, we focus on the case of $\mu\ge2$.
Then it follows from (\ref{ineq_F}) and (\ref{positive_F}) that
\begin{equation}
\label{ineq_F>2}
F''(t)\ge C_4(1+t)^{-(n+\mu/2)(p-1)}|F(t)|^p.
\end{equation}
We shall employ the following lemma now.
\begin{lem}[Takamura\cite{Takamura15}]
\label{lem:odi}Let $p>1,a>0,q>0$ satisfy
\begin{equation}
\label{condi:M}
M:=\frac{p-1}{2}a-\frac{q}{2}+1>0.
\end{equation}
 Assume that $F\in C^2([0,T))$ satisfies
\begin{eqnarray}
& F(t)\ge  At^a &  \mbox{for}\ t\ge T_0,\label{ineq:A}\\
& \d F''(t)\ge B(t+R)^{-q}|F(t)|^p& \mbox{for}\ t\ge0,\label{ineq:B}\\
& F(0)\ge0,\ F'(0)>0,\label{ineq:C}&
\end{eqnarray}
where $A,B,R,T_0$ are positive constants.
Then, there exists a positive constant $C_0=C_0(p,a,q,B)$ such that
\begin{equation}
\label{est:T1}
T<2^{2/M}T_1
\end{equation}
holds provided
\begin{equation}
\label{condi}
T_1:=\max\left\{T_0,\frac{F(0)}{F'(0)},R\right\}\ge C_0 A^{-(p-1)/(2M)}.
\end{equation}
\end{lem}

Due to the lower bound of $F_0$ in (\ref{lower bound}) and the definition of $F(t)$ in \eqref{positive_F},
we have
\begin{equation}
\label{lower bound_F}
F(t)>C_3\e^pt^{2-\mu p+(n-1)(1-p/2)+\mu/2}
\quad\mbox{for}\ t\ge2,
\end{equation}
which is (\ref{ineq:A}) in Lemma \ref{lem:odi} with
\[
A=C_3\e^p,
\quad a=2-\mu p+(n-1)\left(1-\frac{p}{2}\right)+\frac{\mu}{2},
\quad T_0=2.
\]
The inequality (\ref{ineq:B}) with
\[
B=C_4,
\quad q=\left(n+\frac{\mu}{2}\right)(p-1)
\]
follows from (\ref{ineq_F>2}), and
(\ref{ineq:C})  is already established by (\ref{data_Fat0}).
The final step to use Lemma \ref{lem:odi} is to check the sign of $M$. By the assumption that
$p<p_0(n+2\mu)$, we have
\[
M=\frac{\gamma(p,n+2\mu)}{4}>0.
\]

Set
\[
T_0=C_0A^{-(p-1)/(2M)}=C_0C_3^{-2(p-1)/\gamma(p,n+2\mu)}\e^{-2p(p-1)/\gamma(p,n+2\mu)}.
\]
Then, since $F(0)/F'(0)$ is independent of $\e$ by (\ref{data_Fat0}),
one can see that there is an $\e_0=\e_0(f,g,n,p,\mu,R)>0$ such that
\[
T_0\ge\max\left\{2,\frac{F(0)}{F'(0)}\right\}\quad \mbox{for}\ 0<\e\le \e_0 .
\]
This means that $T_1=T_0$ in (\ref{condi}).
Therefore the conclusion of Lemma \ref{lem:odi} implies that the maximal existence
time $T$ of $F(t)$ has to satisfy
\[
T\le C_5\e^{-2p(p-1)/\gamma(p,n+2\mu)}\quad \mbox{for}\ 0<\e \le \e_0,
\]
where
\[
C_5:=2^{8/\gamma(p,n+2\mu)}C_0C_4^{-2(p-1)/\gamma(p,n+2\mu)}>0.
\]
This completes the proof in the case of $\mu\ge2$.
\hfill$\Box$


\section{Proof of Theorem \ref{thm:main} for $\v{0<\mu<2}$}
\par\quad
Before showing the proof of Theorem \ref{thm:main} for $0<\mu<2$, we first prepare the following
lemma:
\begin{lem}
\label{lem:treat-mass}
Suppose that the assumption in Theorem \ref{thm:main} is fulfilled.
Then it holds that
\begin{equation}
\label{F'}
F'(t)>\sqrt{\frac{C_4}{2(p+1)}}(1+t)^{-(n+\mu/2)(p-1)/2}F(t)^{(p+1)/2}
\end{equation}
for $t\ge C_6\e^{-2p(p-1)/\gamma(p,n+2\mu)}$, where we set
\[
C_6:=\left(\frac{\mu(2-\mu)(p+1)}{2C_3^{p-1}C_4}\right)^{1/X}>0
\]
and
\begin{equation}\nonumber
\begin{aligned}
X&:=2-\left(n+\frac{\mu}{2}\right)(p-1)
+(p-1)\left\{2-\mu p+(n-1)\left(1-\frac{p}{2}\right)+\frac{\mu}{2}\right\}\\
&=\frac{\gamma(p,n+2\mu)}{2}.
\end{aligned}
\end{equation}
$C_3$, $C_4$ are the one in (\ref{lower bound}), (\ref{ineq_F}), respectively,
\end{lem}
{\bf Proof.}
Multiplying the both sides of (\ref{ineq_F}) by $(1+t)^2F'(t)>0$ and noting that (\ref{positive_F}),
we get
\[
\begin{array}{l}
\d\frac{(1+t)^2}{2}\left\{\left(F'(t)\right)^2\right\}'
+\frac{\mu(2-\mu)}{8}\{F(t)^2\}'\\
\d\ge C_4(1+t)^{2-(n+\mu/2)(p-1)}F(t)^pF'(t)
\end{array}
\quad\mbox{for}\ t\ge0.
\]
Integration by parts yields that
\[
\begin{array}{l}
\d\frac{(1+t)^2}{2}\left(F'(t)\right)^2
+\frac{\mu(2-\mu)}{8}F(t)^2\\
\d>C_4\int_0^t(1+s)^{2-(n+\mu/2)(p-1)}F(s)^pF'(s)ds
\end{array}
\quad\mbox{for}\ t\ge0.
\]
Noting the assumption on $p$
\[
p\ge p_F(n)>1+\frac{2}{n+\d\frac{\mu}{2}}
\quad\mbox{for}\ \mu>0,
\]
it is easy to get that
\[
2-\left(n+\frac{\mu}{2}\right)(p-1)<0.
\]
And hence we have
\[
\begin{array}{l}
\d\int_0^t(1+s)^{2-(n+\mu/2)(p-1)}F(s)^pF'(s)ds\\
\d\ge(1+t)^{2-(n+\mu/2)(p-1)}\frac{F(t)^{p+1}-F(0)^{p+1}}{p+1}.
\end{array}
\]
Since
\[
p_F(n)=1+\frac{2}{n}>\frac{2}{n+1-\mu}\quad
\mbox{for $n\ge2$ and $0<\mu<2$},
\]
and hence
\[
p>\frac{2}{n+1-\mu}.
\]
This is equivalent to
\[
p\left(2-\mu p+(n-1)\left(1-\frac{p}{2}\right)+\frac{\mu}{2}\right)
>\frac{\gamma(p,n+2\mu)}{2}.
\]
Thus, for $t\ge C_6\e^{-2p(p-1)/\gamma(p,n+2\mu)}\ge 2$ ($\e$ small enough), we have
\begin{equation}
\label{restriction2_new}
C_3\e^pt^{2-\mu p+(n-1)(1-p/2)+\mu/2}\ge2\|f\|_{L^1(\R^n)}\e,
\end{equation}
which implies
\begin{equation}
\label{restriction2}
F(t)\ge 2F(0).
\end{equation}
Hence, it follows from
\[
F(t)^{p+1}-F(0)^{p+1}\ge F(t)^p\left\{F(t)-F(0)\right\}\ge\frac{1}{2}F(t)^{p+1}
\]
that
\begin{equation}
\label{F'-final}
\frac{(1+t)^2}{2}\left(F'(t)\right)^2
+\frac{\mu(2-\mu)}{8}F(t)^2>\frac{C_4(1+t)^{2-(n+\mu/2)(p-1)}F(t)^{p+1}}{2(p+1)}
\end{equation}
for $t\ge C_6\e^{-2p(p-1)/\gamma(p,n+2\mu)}$.

On the other hand, for $t\ge C_6\e^{-2p(p-1)/\gamma(p,n+2\mu)}$, we have
\[
\frac{C_4}{4(p+1)}(1+t)^{2-(n+\mu/2)(p-1)}\{C_3\e^pt^{2-\mu p+(n-1)(1-p/2)+\mu/2}\}^{p-1}
\ge \frac{\mu(2-\mu)}{8}
\]
which gives us
\begin{equation}
\label{restriction1}
\frac{C_4}{4(p+1)}(1+t)^{2-(n+\mu/2)(p-1)}F(t)^{p+1}
\ge\frac{\mu(2-\mu)}{8}F(t)^2
\end{equation}
 by combining \eqref{lower bound_F}.
Therefore, we get (\ref{F'}) from (\ref{F'-final}).

\hfill$\Box$
\par
By (\ref{F'}), it is easy to see that
there is a $\e_0=\e_0(f,g,n,p,\mu,R)>0$ such that
\[
\frac{F'(t)}{F(t)^{1+\delta}}
>\sqrt{\frac{C_4}{2(p+1)}}(1+t)^{-(n+\mu/2)(p-1)/2}F(t)^{(p-1)/2-\delta}
\]
with $0<\delta<(p-1)/2$ holds for
\[
\ t\ge T_1:=C_6\e^{-2p(p-1)/\gamma(p,n+2\mu)}
\quad\mbox{and}\quad 0<\e\le\e_0.
\]
Here we use our lower bound of $F$ in (\ref{lower bound_F}) again to get
\begin{equation}
\label{contradict:F}
\frac{F'(t)}{F(t)^{1+\delta}}
>C_3^{(p-1)/2-\delta}\sqrt{\frac{C_4}{2(p+1)}}\e^{p\{(p-1)/2-\delta\}}t^Y
\quad\mbox{for}\ t\ge T_1,
\end{equation}
where
\[
\begin{array}{ll}
Y
&\d:=\left(\frac{p-1}{2}-\delta\right)
\left\{2-\mu p+(n-1)\left(1-\frac{p}{2}\right)+\frac{\mu}{2}\right\}
-\left(n+\frac{\mu}{2}\right)\frac{p-1}{2}\\
&\d=\frac{\gamma(p,n+2\mu)}{4}-1
-\left\{2-\mu p+(n-1)\left(1-\frac{p}{2}\right)+\frac{\mu}{2}\right\}\delta.
\end{array}
\]
Therefore, taking $\delta$ small enough such that $Y+1>0$, we then have by integrating
(\ref{contradict:F}) over $[T_1,t]$,
\[
\frac{F(T_1)^{-\delta}}{\delta}>\frac{C_3^{(p-1)/2-\delta}}{Y+1}
\sqrt{\frac{C_4}{2(p+1)}}\e^{p\{(p-1)/2-\delta\}}(t^{Y+1}-T_1^{Y+1})
\quad\mbox{for}\ t\ge T_1.
\]
Making use of (\ref{lower bound_F}) with $t=T_1$ in this inequality, we obtain that
\[
1>C_7\e^{p(p-1)/2}T_1^{\gamma(p,n+2\mu)/4-(Y+1)}
(t^{Y+1}-T_1^{Y+1})
\quad\mbox{for}\ t\ge T_1,
\]
where
\[
C_7:=\frac{\delta C_3^{(p-1)/2}}{Y+1}
\sqrt{\frac{C_4}{2(p+1)}}>0
\]
If one sets $t=kT_1$ with $k>1$, then, due to the definition of $T_1$, one has
\[
1>C_7C_6^{\gamma(p,n+2\mu)/4}(k^{Y+1}-1).
\]
Therefore the conclusion of the Theorem \ref{thm:main},
\[
T\le C_8\e^{-2p(p-1)/\gamma(p,n+2\mu)}
\quad\mbox{for}\ 0<\e\le\e_0,
\]
is now established, where
\[
C_8:=\left(1+C_7^{-1}C_6^{-\gamma(p,n+2\mu)/4}\right)^{1/(Y+1)}C_6>0.
\]
This completes the proof in the case of $0<\mu<2$.
\hfill$\Box$

\section*{Acknowledgment}
\par\quad
The first author is partially supported by NSFC(11501273), high level talent project
of Lishui City(2016RC25), the Scientific Research Foundation of the First-Class Discipline of Zhejiang Province(B)(201601), the key laboratory of Zhejiang Province(2016E10007).
The second author is partially supported by the Grant-in-Aid for Scientific Research (C)
(No.15K04964),
Japan Society for the Promotion of Science,
and Special Research Expenses in FY2016, General Topics (No.B21),
Future University Hakodate.
\par
Finally, all the authors are grateful to Prof.M.Reissig
(Technical University Bergakademie Freiberg, Germany)
for his great advice on the classification of the linear problem
for which our first manuscript of arXiv:1701.03232 was not appropriate.


\bibliographystyle{plain}

\end{document}